\documentclass[11pt,runningheads]{llncs}
\usepackage{amsmath}
\usepackage{amsfonts}
\usepackage{amssymb}
\usepackage{latexsym}
\usepackage{epsfig,graphicx}

\newenvironment{Proof}{\noindent\bf{Proof.}\rm}{\hfill$\blacksquare$\bigskip}

\newcommand{\Core}{{\cal{H}}}

\newcommand{\Sat}{{\cal{S}}}

\newcommand{\whp}{\emph{whp}}

\newcommand{\SmallFrac}{e^{-\Omega(n^2p)}}
\newcommand{\SmallFracMN}{e^{-\Omega(m/n)}}

\newcommand{\MainDist}{{\cal{P}}^{{\rm sat}}_{n,m}}
\newcommand{\RandDist}{{\cal{P}}_{n,p}}
\newcommand{\RandDistM}{{\cal{P}}_{n,p}}
\newcommand{\RandDistSAT}{{\cal{P}}^{{\rm sat}}_{n,p}}

\newcommand{\RandDistSATnp}{{\cal{J}}^{{\rm sat}}_{n,p}}

\newcommand{\F}{F}

\newcommand\eps{\varepsilon}

\title{On the random satisfiable process}

\author{
Michael Krivelevich \inst{1}\thanks{Research supported
  in part by a USA-Israel BSF Grant, and by a grant from the
  Israel Science Foundation, and by Pazy Memorial Award.}, Benny Sudakov \inst{2}\thanks{Research supported in part by NSF CAREER award DMS-0546523 and a USA-Israeli BSF grant.}
 and  Dan Vilenchik \inst{3}}

\institute{School of Mathematical Sciences, Tel Aviv University, Tel Aviv, Israel.
\\\email{krivelev@post.tau.ac.il}   \and Department of Mathematics, UCLA, Los Angeles.\\ \email{bsudakov@math.ucla.edu} \and School of
Computer Science, Tel Aviv University, Tel Aviv, Israel.
\\\email{vilenchi@post.tau.ac.il}}

\begin{document}
\pagenumbering{roman}
\maketitle
\begin{abstract}
In this work we suggest a new model for generating random
satisfiable $k$-CNF formulas. To generate such formulas -- randomly
permute all $2^k\binom{n}{k}$ possible clauses over the variables
$x_1,\ldots,x_n$, and starting from the empty formula, go over the
clauses one by one, including each new clause as you go along if
after its addition the formula remains satisfiable. We study the
evolution of this process, namely the distribution over formulas
obtained after scanning through the first $m$ clauses (in the random
permutation's order).

Random processes with conditioning on a certain property being respected are widely
studied in the context of graph properties. This study was
pioneered by
Ruci\'nski and Wormald in 1992 for graphs with a fixed degree sequence,
and also by Erd\H{o}s, Suen, and Winkler in 1995 for triangle-free
and bipartite graphs. Since then many other graph properties were studied such as planarity and $H$-freeness.
Thus our model is a natural extension of this approach to the satisfiability setting.

Our main contribution is as follows. For $m \geq cn$, $c=c(k)$ a sufficiently large
constant, we are able to characterize the structure of the solution
space of a typical formula in this distribution. Specifically, we
show that typically all satisfying assignments are essentially clustered in
one cluster, and all but $e^{-\Omega(m/n)} n$ of the variables take
the same value in all satisfying assignments. We also describe a
polynomial time algorithm that finds $\whp$ a satisfying assignment
for such formulas. %Our result is in agreement with what is already
%known about the structure of the solution space of formulas drawn
%from other distributions over satisfiable $k$-CNF formulas (for the
%same order of clause-variable ratio) such as the uniform
%distribution and the planted one.

\end{abstract}
\newpage
\pagenumbering{arabic}
\section{Introduction}
Constraint satisfaction problems play an important role in many
areas of computer science, e.g. computational complexity theory
\cite{Cook71}, coding theory \cite{Gallager}, and artificial
intelligence \cite{Pearl}, to mention just a few. The main
challenge is to devise efficient algorithms for finding satisfying
assignments (when such exist), or conversely to provide a certificate
of unsatisfiability. One of the best known examples of a
constraint satisfaction problem is $k$-SAT, which is the first to
be proven as NP-complete. Although satisfactory approximation
algorithms are known for several NP-hard problems, the problem of
finding a satisfying assignment (if such exists) is not amongst
them. In fact, H{\aa}stad \cite{Hastad01} proved that it is
NP-hard to approximate MAX-3SAT (the problem of finding an
assignment that satisfies as many clauses as possible) within a
ratio better than 7/8.

\medskip

In trying to understand the inherent hardness of the problem, many
researchers analyzed structural properties of formulas drawn from
different distributions. One such distribution is the \textbf{uniform} distribution where instances
are generated by picking $m$ clauses uniformly at random out of all $2^k\binom{n}{k}$ possible clauses. Although
many problems still remain unsolved, in general this distribution seems to be quite well understood (at least for some values of $m$ and $k$). This is also true for the
\textbf{planted $k$-SAT} model, where one first fixes some assignment $\psi$ to the variables and then picks $m$ clauses uniformly at random out of all $(2^k-1)\binom{n}{k}$ clauses
satisfied by $\psi$. Comparatively, much less is known for variants of these distributions where extra conditions are imposed. These conditions
distort the randomness in such a way that the ``standard" methods and tools employed to analyze the original distributions are a-priori of little use in the new setting. Our
work concerns the latter.

\subsection{Our Contribution}\label{sec:ourCont}
In this work we suggest a new model for generating random
satisfiable $k$-CNF formulas. To generate such formulas -- randomly
permute all $2^k\binom{n}{k}$ possible clauses over the variables
$x_1,\ldots,x_n$, and starting from the empty formula, go over the
clauses one by one, including each new clause as you go along if
after its addition to the formula, the formula remains satisfiable.
We study the evolution of this process, namely the distribution over
formulas obtained after scanning through the first $m$ clauses (in
the random permutation's order); we use $\MainDist$ to denote this
distribution. Clearly, for every $m$, all formulas in $\MainDist$
are satisfiable (as every clause is included only if the so-far
obtained formula remains satisfiable).

Random processes with conditioning on a certain property being respected are widely
studied in the context of graph properties. This study was pioneered by
Ruci\'nski and Wormald in 1992 \cite{RucinskiW92} for graphs
with a fixed degree sequence, and also by Erd\H{o}s, Suen, and Winkler in
1995 for triangle-free
and bipartite graphs \cite{ErdosMaxGraph}. Since then many other graph properties were studied such as planarity \cite{PlanarGraphProcess}, $H$-freeness \cite{HFreeGraphs} and also the property of being intersecting in the context of hypergraphs~\cite{IntersectingHyper}.
Thus our model is a natural extension of this approach to the
satisfiability setting. The main difficulty when dealing with
these restricted processes is that the edges of the random graph (and the
clauses of the random $k$-CNF formula) are no longer independent due to
conditioning. Thus the rich methods that have been developed to understand
the ``classical" random graph models, $G_{n,p}$ for example, do not carry
over, at least not immediately, to the restricted setting.

Quite frequently in restricted random processes, the typical size
of a final graph or formula (after all edges/clauses have been
scanned) is a fascinating subject of study.
This is however {\em not} the case here, as it is quite easy to
see that deterministically the final random formula will have
$(2^k-1){n\choose k}$ clauses and a unique satisfying assignment.
Therefore, in the setting under consideration here the process
itself (i.e. a typical development of a restricted random formula
and of its set of satisfying assignments as the number of scanned clauses
$m$ grows) is much more interesting than the final result, and
indeed in this paper we will study the development of a random
satisfiable formula.

As it turns out,
if $m$ is chosen so that almost all $k$-CNF formulas with $m$
clauses over $n$ variables are satisfiable, then $\MainDist$ is
statistically close to the uniform distribution over such formulas
since $\whp$ none of the $m$ clauses will be rejected (writing
$\whp$ we mean with probability tending to $1$ as $n$ goes to
infinity). Therefore if this is the case, then the clauses are practically
independent of each other, and the ``usual" techniques apply. Remarkable phenomena occurring in the uniform
distribution are {\bf phase transitions}. With respect to the
property of being satisfiable, such a phase transition takes place
too. More precisely, there exists a threshold $d=d(n,k)$ such that
almost all $k$-CNF formulas over $n$ variables with $m$ clauses such
that $m/n>d$ are not satisfiable, and almost all formulas with
$m/n<d$ are~\cite{Friedgut}. Thus, while $\MainDist$ is
statistically close to the uniform distribution for $m/n$ below the
threshold, it is not clear how does a typical $\MainDist$ instance
look like when crossing this threshold (which is conjectured to be
roughly 4.26 for 3SAT), and whether there exists a polynomial time
algorithm for finding a satisfying assignment for such instances.
%One main difficulty in studying $\MainDist$ in this regime is the
%fact that the ``standard" techniques which are used to study random
%distributions where clauses (edges, in the graph setting) are chosen
%independently (or nearly so) no longer carry over to $\MainDist$,
%$m/n$ above the satisfiability threshold.

In this work we analyze $\MainDist$ when $m/n$ is some sufficiently
large constant \emph{above} the satisfiability threshold.
The first part of our result is characterizing the structure of the solution space of a typical
formula in $\MainDist$. By the ``solution space" of a formula we mean the set of all satisfying assignments (which is a subset
of all $2^n$ possible assignments). Formally,

\begin{theorem}\label{thm:StructOfSolutionSpace} Let $\F$ be random $k$-CNF from $\MainDist$,
$m/n \geq c$, $c=c(k)$ a sufficiently large constant. Then $\whp$
$\F$ enjoys the following properties:
\begin{enumerate}
  \item All but $\SmallFracMN n$ variables are frozen.
  \item The formula induced by the non-frozen variables decomposes into connected components of at most logarithmic size.
  \item Letting $\beta(\F)$ be the number of satisfying assignments of $\F$,
    we have $\frac{1}{n}\log\beta(F)=\SmallFracMN$.
\end{enumerate}
\end{theorem}
By a {\em frozen variable} we mean a variable that takes the same value in all satisfying assignments.
Notice that the third item in Theorem \ref{thm:StructOfSolutionSpace} follows directly from the first.
One immediate corollary of this theorem is:
\begin{corollary} Let $\F$ be random $k$-CNF from $\MainDist$,
$m/n \geq c\log n$, $c=c(k)$ a sufficiently large constant. Then
$\whp$ $\F$ has only one satisfying assignment.
\end{corollary}
The corollary follows from the third item in Theorem
\ref{thm:StructOfSolutionSpace} since $\SmallFracMN=o(n^{-1})$
for $m/n\geq c\log n$, and therefore $\log \beta(F)=o(1)$, or in
turn, $\beta(F)=1+o(1)$.

\medskip

The characterization given by Theorem
\ref{thm:StructOfSolutionSpace} is in sharp contrast with the structure
of the solution space of $\MainDist$ formulas with $m/n$ just below the threshold. Specifically, the
\emph{conjectured} picture, some supporting evidence of which was
proved rigorously for $k\geq 8$
\cite{AchiRicciTers06,ClusteringPhysicists,AminAchi}, is that
typically random $k$-CNF formulas in the near-threshold regime have
an exponential number of {\bf clusters} of satisfying assignments.
While any two assignments in distinct clusters disagree on at least
$\eps n$ variables, any two assignments within one cluster coincide
on $(1-\eps)n$ variables. Furthermore, each cluster has a linear
number of frozen variables (frozen w.r.t. all satisfying assignments within that
cluster). This structure seems to make life hard for most known SAT heuristics.
One explantation seems to be that
the algorithms do not ``steer'' into one cluster but rather try to
find a ``compromise'' between the satisfying assignments in distinct
clusters, which actually is impossible.

Complementing this picture \emph{rigorously}, we show that a typical
formula in $\MainDist$ (in the above-threshold regime) can be solved
efficiently. Formally,

\begin{theorem}\label{thm:PlytimeAlg} There exists a deterministic polynomial time algorithm that
$\whp$ finds a satisfying assignment for $k$-CNF formulas from
$\MainDist$, $m/n\geq c$, $c=c(k)$ a sufficiently large constant.
\end{theorem}
Our proof of Theorem \ref{thm:PlytimeAlg} is constructive in the sense that we explicitly
describe the algorithm.

\begin{remark}\label{rem:depOnk}
Observe that in both theorems we have $m/n \geq c(k)$, $c$ some function of $k$. We assume that $k$
is fixed, and therefore $c(k)$ is some constant. The true dependency is given by $c(k)=c_02^k$ where $c_0$ is
some moderate universal constant, say 100. The exponential dependency on $k$ is somewhat inevitable as the
satisfiability threshold itself scales exponentially with $k$ (asymptotically $2^k\ln 2$). In this work we do not
go to such fine details as determining the constant $c_0$, though
this task is perhaps manageable.
\end{remark}

\begin{remark}\label{rem:k-col}
Another natural problem to study is $k$-colorability. Similar to
random $k$-CNF formulas, the random graph $G_{n,p}$ also goes
through a phase transition w.r.t. the property of being
$k$-colorable, as $np$ grows. Analogously to the random $k$-CNF
process that we defined, one can consider a restricted random graph
process.
Specifically, randomly order all $\binom{n}{2}$ edges of the graph,
go over them in that order and include each new edge as long as the
resulting graph remains $k$-colorable. Some of the results that we
have for $k$-SAT extend to the $k$-colorability process. A
more thorough discussion is given in Section \ref{sec:k-col}.
\end{remark}
%\medskip
%
%Combining Theorems \ref{thm:StructOfSolutionSpace} and \ref{thm:PlytimeAlg} suggests an interesting connection between the
%hardness (even experimental one) of a certain clause-density regime
%and the typical structure of the solution space at
%that point. Specifically, our results show that when at the
%single-cluster regime the problem is ``easy", in contrast to the
%``hard" multi-clustered near-threshold regime.
%
%\medskip
\subsection{Related Work and Techniques}
Almost all polynomial-time heuristics suggested so far for
random instances (either SAT or graph optimization problems) were
analyzed when the input is sampled according to a planted-solution
distribution, or various semi-random variants thereof. Alon and
Kahale \cite{AlonKahale97} suggest a polynomial time algorithm
based on spectral techniques that $\whp$ properly $k$-colors a
random graph from the planted $k$-coloring distribution (the
distribution of graphs generated by partitioning the $n$ vertices
into $k$ equally-sized color classes, and including every edge
connecting two different color classes with probability $p=p(n)$),
for graphs with average degree greater than some constant. In the
SAT context, Flaxman's algorithm, drawing on ideas from \cite{AlonKahale97}, solves $\whp$ planted 3SAT instances
where the clause-variable ratio is greater than some constant. Also \cite{TechReport,WP,ExpectedPoly3SAT} address the planted 3SAT distribution.

\medskip

On the other hand, very little work was done on non-planted
distributions, such as $\MainDist$. In this context one can mention
a work of Chen \cite{Chen03} who provides an \emph{exponential}
time algorithm for the uniform distribution over satisfiable $k$-CNF formulas with exactly $m$ clauses where $m/n$
is greater than some constant. Ben-Sasson et al. \cite{EBSPlanted} also study this
distribution but with $m/n=\Omega(\log n)$, a regime where
the uniform distribution  and the planted distribution essentially coincide
(since typically there is only one satisfying assignment), and leave
as an open question whether one can characterize the regime
$m/n=o(\log n)$. This question was resolved in
\cite{UniformSAT} (and in \cite{UniformCol} for the uniform
distribution over $k$-colorable graphs).

While some of the ideas suggested in these works have proven to be
instrumental for our setting, most
of their analytical methods break when considering $\MainDist$. In
$\MainDist$ not only do clauses depend on each other (unlike the
planted distribution where clauses are chosen independently), but
the order in which they are introduced also plays a role (which is
not the case in the uniform distribution studied in
\cite{UniformSAT}, although the clauses are not chosen
independently). Therefore we had to come up with new analytical
tools that might be of interest in other settings as well.

\subsection{Paper's Structure}

The rest of the paper is structured as follows. In Section
\ref{sec:PropOfRandomInst} we discuss relevant structural properties
that a typical formula in $\MainDist$ possesses, the proofs of some properties are postponed to Sections \ref{sec:ExpndrProof} and \ref{sec:ProofOfPropSizeOfConnectedComp}. One consequence of this discussion will be a proof of Theorem
\ref{thm:StructOfSolutionSpace}. We then prove Theorem
\ref{thm:PlytimeAlg} in Section \ref{sec:ProofOfThm} by presenting
an algorithm and showing that it meets the requirements of Theorem
\ref{thm:PlytimeAlg}. In Section \ref{sec:k-col} we discuss the $k$-colorability setting
(mentioned in Remark \ref{rem:k-col}) more elaborately, and concluding remarks are given in Section
\ref{sec:Discussion}.

\medskip

To simplify the presentation we shall address, in what follows, only
the case $k=3$. The case of general $k$ easily follows from the same
arguments (taking $m/n \geq c(k)$, $c(k)$ as mentioned in Remark \ref{rem:depOnk}).

\section{Properties of a Random $\MainDist$ Instance}\label{sec:PropOfRandomInst}
This section contains the technical part of the paper.
In it we analyze the structure of a typical formula in
$\MainDist$. Here and throughout we think of $m$ as $cn$, $c$ at least some sufficiently large constant.

\subsection{Preliminaries and Techniques}
When analyzing some structural properties of a random instance in
$\MainDist$ it will be more convenient to analyze the same property
under a somewhat different distribution, and then to go back to
$\MainDist$ (maybe pay some factor in the estimate).

The variation we consider is $\RandDistSAT$ and is defined as follows:
permute at random all possible $M=8\binom{n}{3}$ clauses, go over the
clauses in the permutation's order and include each clause with probability
$p=m/M$ if also its addition leaves the instance satisfiable. Let
$\RandDist$ be defined similarly, just without the conditioning
(i.e., all clauses chosen at random are included in the formula,
thus making it not necessarily satisfiable).

\begin{lemma} $\MainDist=\RandDistSAT|\{\text{ exactly $m$ clauses
were chosen}\}$.
\end{lemma}
\begin{Proof} To generate $\F$ in $\MainDist$ one first picks a
random permutation of the clauses and then scans one by one the
first $m$ clauses, skipping clauses whose addition will make
the instance unsatisfiable. The key point is to notice that any
ordered $m$-tuple of clauses is equally likely to be chosen as the
first $m$ clauses. This is exactly the case in $\RandDistSAT$ when
conditioning on the fact that exactly $m$ clauses were chosen --
any set of $m$ clauses is equally likely, and also any
permutation of them.
\end{Proof}

\begin{lemma}\label{lem:Main2RandSATnp} Set $M=8\binom{n}{3}$. For any property $A$, if $p=m/M$ then $Pr^{\MainDist}[A] \leq
O(\sqrt{m})\cdot Pr^{\RandDistSAT}[A]$.
\end{lemma}
\begin{Proof}
Let $X$ be a random variable counting the number of clauses whose coin toss was successful. $X$ is distributed
$Binom(8 \binom{n}{3},p)$, and therefore $E[X]=m$. Standard calculations show that $Pr[X = m] = \Omega(m^{-0.5})$.

$$Pr^{\MainDist}[A]=Pr^{\RandDistSATnp}[A|X=m]=\frac{Pr^{\RandDistSATnp}[A \wedge X=m]}{Pr^{\RandDistSATnp}[X=m]} \leq O(\sqrt{m})\cdot Pr^{\RandDistSAT}[A].$$
\end{Proof}

\begin{remark}\label{rem::RedfinigWhp}
In the remainder of the section we analyze $\RandDistSAT$ instead of $\MainDist$. When we use the expression ``with high probability" (abbreviated $\whp$) we will always mean with probability $1-o(m^{-1/2})$. Lemma \ref{lem:Main2RandSATnp} will then imply that we can switch back to $\MainDist$ and still the property holds with probability $1-o(1)$. We will actually prove that all the properties hold with probability $1-o(n^{-3})$ which is always at least $1-o(m^{-1/2})$ since $m=O(n^3)$.
\end{remark}

\subsection{The Discrepancy Property}\label{sec:discrepacny}
A well known result in the theory of random graphs is that a random
graph $\whp$ will not contain a small yet unexpectedly dense
subgraph. This is also the case for $\RandDist$ (when considering
the graph induced by the formula). In general, discrepancy
properties play a fundamental role in the proof of many important
structural properties such as expansion, the spectra of the
adjacency matrix, etc., and indeed in our case the discrepancy property
plays a major role both in the algorithmic perspective and in the
analysis of the clustering phenomenon. The following discussion
rigorously establishes the above stated fact.

\begin{definition}\label{def:proportional} We say that a 3CNF
formula $\F$ on $n$ variables is \textbf{$\rho$-proportional} if there exists \emph{no} set $U$ of
variables such that:
\begin{itemize}
    \item $|U|\leq n/10^6$,
    \item There are at least $\rho\cdot |U|$ clauses in $\F$ each containing at least two variables
    from $U$.
\end{itemize}
\end{definition}
(We say that a clause $C$ contains a variable $x$ if $x$ appears in $C$
either as $x$ or as $\bar{x}$, in this context
we do not differentiate between the two cases).

\begin{proposition}\label{prop:NoDenseSubgraphs} Let $\F$ be distributed according to $\RandDist$ with
$n^2p\geq d$, $d$ a sufficiently large constant, and set
$\rho=n^2p/5500$. Then $\whp$ $\F$ is $\rho$-proportional.
\end{proposition}
\begin{remark} To see how Proposition \ref{prop:NoDenseSubgraphs}
corresponds to the random graph context, consider the graph induced
by the formula $\F$ (the vertices are the variables, and two
variables are connected by an edge if there exists some clause containing them
both) and observe that every clause that contains at least two
variables from $U$ contributes an edge to the subgraph induced by
$U$. Thus if we have many such clauses, this subgraph will be
prohibitively dense. Since $\F$ is random so is its induced graph,
and therefore the latter will typically not occur.
\end{remark}
\begin{Proof}
The probability that a random formula $\F$ in $\RandDist$ contains a set $U$ of variables of size $u$ that violates
proportionality is at most (using the union bound):
$$\sum_{u=1}^{n/10^6}
\binom{n}{u}\cdot\binom{8n\binom{u}{2}}{un^2p/5500}\cdot
p^{un^2p/5500}=o(n^{-3}).$$
The first term accounts for the possible ways of choosing the variables of $U$, the second is to choose
the $un^2p/5500$ clauses that contain at least two variables from $U$ (out of at most $8n\binom{u}{2}$ possible
ones), and the last term is just the probability of the chosen
clauses to actually appear in $F$.
To bound this sum we use the fact that $u \leq n/10^6$, the fact that $n^2p$ can be arbitrarily large (constant), and the following standard estimate for the binomial coefficient:
$$\binom{n}{x}\leq \left(\frac{en}{x}\right)^x.$$
\end{Proof}
\begin{corollary}\label{cor:NoDensSub} Let $\F^*$ be distributed according to $\RandDistSAT$ with
$n^2p\geq d$, $d$ a sufficiently large constant. Then $\whp$
$\F^*$ is $n^2p/5500$-proportional.
\end{corollary}
The corollary follows easily by observing that the proportionality
property is monotonically decreasing.

\subsection{Crude Characterization of the Solution Space's Structure}
In this section we make the first step towards proving Theorem
\ref{thm:StructOfSolutionSpace} (clustering). We give a rather crude
characterization of the structure of the solution space of a typical
instance in $\RandDistSAT$. This characterization will be refined in
the sequel.

\begin{definition}\label{def:concentrated} A 3CNF $\F$ is called $r$-\textbf{concentrated} if every
two satisfying assignments $\psi_1,\psi_2$ of $\F$ are at Hamming
distance at most $r$ from each other.
\end{definition}

\begin{proposition} \label{prop:Concentration}  Let $\F^*$ be distributed according to $\RandDistSAT$ with
$n^2p \geq d$, $d$ a sufficiently large constant, let $\rho=30/(n^2p)$ then $\whp$ $\F^*$ is
$\rho n$-concentrated.
\end{proposition}
An immediate corollary of this proposition is that typically all
satisfying assignments of a $\RandDistSAT$ instance can be enclosed
in a ball of radius $30/(np)$ in $\{0,1\}^n$. This gives a
``first-order" characterization of the structure of the solution
space.

\medskip

\begin{Proof} Fix two assignments $\varphi$ and $\psi$ at distance
$\alpha n$, and let us bound $Pr[\varphi \text{ and } \psi \text{
satisfy } F^*]$. Assume w.l.o.g. that, say, $\varphi$ is the
all-TRUE assignment. We shall now upper bound the probability of a set of clauses in $\RandDist$
that may result in an instance $F^*$ that is satisfied by both assignments.
In particular a clause
of the form $C_1=(x \vee \bar{y} \vee \bar{z})$, where $x$ is a
variable on which $\varphi$ and $\psi$ disagree, and $y,z$ are variables on which both agree, cannot be chosen to $\RandDist$. Let us call such a clause a type 1 clause. If a type 1 clause appears is included, then either it is included in $F^*$, and then $\psi$ cannot be a satisfying
assignment, or it is rejected and then $\varphi$ is already at this point not a satisfying assignment. The same applies for clauses of the form $C_2=(s \vee w \vee t)$, where on all three variables, $s,w,t$, both assignments disagree -- call them type 2.
It remains to upper bound the probability of a $\RandDist$ instance that does not contain type 1 and type 2 clauses. There are $\alpha n\binom{(1-\alpha)n}{2}$ type 1 clauses and $\binom{\alpha n}{3}$ type 2 clauses.
The probability of none being chosen is
\begin{equation}\label{eq:ConcentrationEq}
(1-p)^{\alpha n\binom{(1-\alpha)n}{2}+ \binom{\alpha n}{3}} \leq \exp\{-p \cdot \left(\alpha n\binom{(1-\alpha)n}{2}+ \binom{\alpha n}{3}\right)\}.
\end{equation}
If $30/n^2p \leq \alpha \leq 1/2$ then
$$p \cdot \alpha n\binom{(1-\alpha)n}{2} \geq p \alpha n \cdot
n^2/8 \geq 3n.$$
If $\alpha \geq 1/2$ then

$$p \cdot \binom{\alpha n}{3} \geq n \cdot n^2p/48 \geq 3n.$$
In the last inequality we use the fact that $n^2p$ can be
arbitrarily large (specifically, greater than 144).
In any case, the expression in (\ref{eq:ConcentrationEq}) is at most
$5^{-n}$. Since we have no more than $4^n$ ways of choosing the pair
$\varphi,\psi$, we deduce using the union bound that $\whp$ no
such ``bad" pair exists.
\end{Proof}

\subsection{The Core Variables}\label{sec:CoreVertices}
We describe a subset of the variables, referred to as the
\emph{core variables}, which plays a crucial role in the
understanding of $\RandDistSAT$. A variable is said to be frozen
in $\F$ if in every satisfying assignment it takes the same
value. The notion of a core captures this phenomenon. In
addition, a core typically contains all but a small (though
constant) fraction of the variables. This implies that a large
fraction of the variables is frozen, a fact which must leave
imprints on various structural properties of the formula. These
imprints allow efficient heuristics to recover a satisfying
assignment of the core. A second implication of this is an upper
bound on the number of possible satisfying assignments, and on the
distance between every such two. Thus the notion of a core plays a
key role in obtaining a characterization of the cluster structure
of the solution space.

\medskip

%A somewhat different notion of core is also defined
%in~\cite{AlonKahale97,flaxman,UniformSAT,UniformCol}. However, as we
%already mentioned, our distribution significantly differs from the
%ones discussed in the latter papers. In particular, the techniques
%developed in these papers to prove the existence of a large core do
%not carry over to our setting (at least not in any way that we are
%aware of).

Let us now proceed with a rigorous definition of a core. Before
doing so, we take a long detour on expanding sets.

\begin{definition}(support)
Given a 3CNF formula $\F$ and some assignment $\psi$ to the
variables, we say that a variable $x$ \textbf{supports} a clause $C$
(in which it appears) w.r.t. $\psi$ if $x$ is the only variable whose literal
evaluates to true in $C$ under $\psi$.
\end{definition}

\begin{definition}\label{def:expanding}(expanding set) Given a
3CNF formula $\F$
%on $n$ variables with $m$ clauses
and an assignment $\psi$ to the variables (not
necessarily satisfying), a set of variables $Z$ is called
\textbf{$t$-expanding} in $\F$ w.r.t. $\psi$ if every variable $x\in Z$
supports at least $t$ clauses in $\F[Z]$ w.r.t. $\psi$.
\end{definition}
$\F[Z]$ stands for the subformula of $\F$ containing the clauses
where all three variables belong to $Z$. The following proposition
illustrates the usefulness of Definition \ref{def:expanding}.
\begin{proposition}\label{prop:UniqueOfSat}
Let $\F$ be a 3CNF formula on $n$ variables
%with $m$ clauses
and let $Z$ be a $t$-expanding set w.r.t. some assignment $\psi$.
If in addition:
\begin{itemize}
    \item $\psi$ satisfies $\F$,
    \item $\F$ is $n/10^6$-concentrated (Definition
    \ref{def:concentrated}),
    \item $\F$ is $t$-proportional (Definition \ref{def:proportional}),
\end{itemize}
then the variables in $Z$ are frozen in $\F$.
\end{proposition}
\begin{Proof} By contradiction, let $\psi$ be the satisfying assignment
w.r.t. which $Z$ is defined and let $\psi'$ be some satisfying
assignment of $\F$ such that there exists a non-empty set $U \subseteq Z$ of variables for
which $\forall x\in U,\psi(x)\neq\psi'(x)$ (if for every $\psi'$
it holds that $U=\emptyset$ then we are done). Take $x\in U$ and
consider all the clauses that $x$ supports w.r.t. $\psi$ in
$\F[Z]$. It must be that every such clause contains at least another
variable $y$ on which $\psi$ and $\psi'$ disagree (since every
such clause is satisfied by $\psi'$ but the literal corresponding to $x$ is
false under $\psi'$). Therefore $y$ belongs to $U$ by definition. We conclude that
there exists a set $U$ of variables and $t\cdot |U|$
clauses each containing at least two variables from $U$ (no clause
was counted twice since the supporter of a clause is unique by
definition). Further, we assumed that $\F$ is $n/10^6$-concentrated
and therefore $|U|\leq n/10^6$. Combining the latter two facts we
derive a contradiction to the $t$-proportionality of $\F$.
\end{Proof}
\begin{proposition}\label{prop:ExpanderSize} Let $\F$ be distributed according to $\RandDistM$ with
$n^2p\geq d$, $d$ a sufficiently large constant. Then $\whp$ there
exists an integer $t=t(n,p)>0$, a set $Z$ of variables, and an assignment $\psi$ such that:
\begin{itemize}
    \item $Z$ is $t$-expanding w.r.t. $\psi$,
    \item $|Z|=(1-\SmallFrac)n$.
    \item $\F$ is $t/10$-proportional,
    \item $\psi$ satisfies $\F^*$,
    \item $\F^*$ is $n/10^6$-concentrated,
\end{itemize}
\end{proposition}
The complete proof of this proposition is deferred to Section
\ref{sec:ExpndrProof}.
\begin{corollary}\label{cor:FrozenInSatFormula} The set
$Z$ promised in Proposition \ref{prop:ExpanderSize} is frozen in $\F^*$, and furthermore $Z$ is $t$-expanding w.r.t.
every satisfying assignment of $\F^*$.
\end{corollary}
\begin{Proof}
To see why $Z$ is frozen, let $S$ be the set of clauses in $\F$ that
are supported w.r.t. $\psi$. First observe that $\F^*$ is
$t$-proportional as it is a subformula of $\F$ (and $\F$ is $t/10$-proportional
and therefore also $t$-proportional). Furthermore $S$ is
contained in $\F^*$. This is because $\psi$ is a satisfying
assignment of $\F^*$ throughout the entire generating process, thus
every clause in $S$ that arrives is not rejected. Therefore $Z$ is also
$t$-expanding in $F^*$  w.r.t. $\psi$.
Finally apply Proposition \ref{prop:UniqueOfSat} to $\F^*$.
The second part of the corollary is immediate from the fact that $Z$
is frozen.
\end{Proof}
\begin{definition}\label{def:SelfContained}(self-contained sets) Given a
3CNF formula $\F$
%on $n$ variables with $m$ clauses
we say that a set of variables $Z$ is
\textbf{$r$-self-contained} in $\F$ if every variable $x\in Z$
appears in at most $r$ clauses in $\F\setminus\F[Z]$.
\end{definition}
Finally, we are ready to define a core.
\begin{definition}\label{def:core}(core)
A set of variables $\Core$ is called a $t$-\textbf{core} of $\F$
w.r.t. an assignment $\psi$ if $\Core$ is $t$-expanding in $F$
w.r.t $\psi$ and also $(t/3)$-self-contained in $\F$.
\end{definition}
The property of being self-contained is necessary for the algorithmic part (the proof of Theorem \ref{thm:PlytimeAlg}, at least
as our analysis proceeds).
\begin{proposition}\label{prop:CoreSize} Let $\F^*$ be distributed according to $\RandDistSAT$ with
$n^2p\geq d$, $d$ a sufficiently large constant. Then $\whp$ there
exists an integer $t=t(n,p)>0$, a satisfying assignment $\psi$ of $F^*$, and a $t$-core $\Core$ w.r.t. $\psi$ such that:
\begin{itemize}
  \item $|\Core|=(1-\SmallFrac)n$,
  \item $\Core$ is frozen in $\F^*$,
  \item $\F^*$ is $t/10$-proportional.
\end{itemize}
\end{proposition}
The proof of this proposition is best understood in the context of
the proof of Proposition \ref{prop:ExpanderSize}. Therefore the
proof appears in Section \ref{sec:CoreSizeProof}.

%\begin{corollary}\label{cor:CoreSizeInSATDist} Let $\F^*$ be distributed according to $\RandDistSAT$ with
%$m\geq C_0n$, $C_0$ a sufficiently large constant. Then $\whp$ there exists a $t$-core $\Core$ enjoying the same
%properties stated in the proposition.
%\end{corollary}
%\begin{Proof}(corollary) Let $\Core$ be the $t$-core promised by Proposition \ref{prop:CoreSize}. Observe that since $\Core$ is $t$-expanding
%w.r.t. a satisfying assignment of $\F^*$, then all the supported
%clauses are also included in $\F^*$. Therefore $\Core$ is also
%$t$-expanding in $\F^*$. As for being $t$-self-contained -- this
%property is closed under removal of clauses. Therefore it must hold
%for $\F^*$ as well.
%\end{Proof}
\begin{remark}\label{rem:MaximalityOfCore}
Observe that if there exist two $t$-cores $\Core_1$ and $\Core_2$
that satisfy the conditions of Proposition \ref{prop:CoreSize}, then
also their union $\Core_1 \cup \Core_2$ is a $t$-core (since the
core variables are frozen). Therefore we may speak of a unique
maximal $t$-core. From now on, when we refer to a $t$-core, we mean
the maximal one. Note that this maximal core is also frozen by
Proposition \ref{prop:UniqueOfSat}. Therefore it can serve as a $t$-core
for {\em any} satisfying assignment of $F$ and thus is effectively uniquely
defined by the formula.
\end{remark}

\subsection{Satellite Variables}\label{sec:Sattelite}
In this section we isolate another set of variables which we call satellite variables. As it turns out, to prove Theorems \ref{thm:StructOfSolutionSpace} and \ref{thm:PlytimeAlg}, it is enough to distinguish between core and satellite variables
and all other variables in $V$. Let us start with a formal definition of a satellite variable.

\begin{definition}\label{def:Satellite} Given a formula $F$ with a core set $\Core$ w.r.t. to an assignment $\varphi$, a variable $x$ is called a \emph{$0$-satellite} with respect to $\Core$ if $x \in \Core$. A variable $x$ is called an \emph{$i$-satellite} if $F$ contains a clause of the form $(x \vee \ell_{z_1} \vee \ell_{z_2})$ or $(\bar{x} \vee \ell_{z_3} \vee \ell_{z_4})$
where for every $j=1,2,3,4$, $z_j$ is a $b$-satellite for $b<i$, and $\varphi(\ell_{z_j})=FALSE$, moreover at least one of $z_j$ is an $(i-1)$-satellite. We say that $x$ is a satellite variable if it is $b$-satellite for some number $b \geq 1$.
\end{definition}
In this definition, $\ell_z$ stands for a literal corresponding to
a variable $z$ (i.e. $\ell=z$ or $\ell=\bar{z}$). Observe that
if $\Core$ is frozen in $F$ then $\Core \cup \Sat$ is frozen as
well (this follows from a simple inductive argument).

\medskip

Before we formally state the property involving the satellite variables we introduce some additional notation. The connected components of a formula $\F$ are the sub-formulas $\F[C_1],\ldots,\F[C_k]$,
where $C_1,C_2,\ldots,C_k$ are the connected components in the
graph $G_\F$ induced by $\F$ (the vertices of $G_\F$ are the variables, and two
variables are connected by an edge if there exists some clause containing them
both). Given a set of variables  $A$ and an assignment $\varphi$ we denote by
$\F_{out}(A,\varphi)$ the subformula of $\F$ which is the
outcome of the following procedure: set the variables in $A$ according to $\varphi$ and simplify $\F$ (by simplify we mean remove every clause that contains a TRUE literal, and remove FALSE literals from the other clauses).

\begin{proposition}\label{prop:SizeOfConnectedComp} Let $\F^*$ be distributed according to $\RandDistSAT$ with
$n^2p\geq d$, $d$ a sufficiently large constant. There $\whp$ exists an integer $t=t(n,p)>0$, a satisfying assignment $\psi$ of $F^*$, and a $t$-core $\Core$ w.r.t. $\psi$ such that:
\begin{itemize}
  \item $|\Core| \geq (1-\SmallFrac)n$.
  \item $\F^*$ is $t/10$-proportional.
  \item Let $\Sat$ be its satellite variables, $\Core \cup \Sat$ are frozen in $\F^*$,
  \item The largest connected component in
  $\F^*_{out}(\Core \cup \Sat,\psi)$ is of size at most $\log n$.
\end{itemize}
\end{proposition}

The new addition compared with Proposition \ref{prop:CoreSize} is the fact that we characterize the structure of the formula induced by the variables not in $\Core \cup \Sat$.

Our proof strategy is the following. Expose the first part of the
random formula $F$ and consider a $t$-core $\Core$ promised $\whp$
by Proposition \ref{prop:CoreSize}. We look at a ``large" connected
component outside the core (if none exists then we are done) and consider
the following ``shattering" procedure. Expose the second part of
the random formula, and suppose for the time being that the core does not
change (even if new clauses are included in $F^*$). Let $x$ be a non-core
variable after the first part, which lies in a spanning tree of a large
connected component. The key observation is that when resuming the random
clause process, $x$ becomes a satellite variable with high (constant)
probability, in which case the spanning tree splits into parts. Since the
tree is large, it contains many variables $x$, and therefore
with very high probability at least one of them will become a satellite
variable and shatter the tree. Finally, it remains to upper bound the
number of possible large trees vs. the probability that such a
tree does not survive. The complete proof is given is
Section \ref{sec:ProofOfPropSizeOfConnectedComp}.

One problem with the approach we just described is that we assumed that the core $\Core$ established after the first round does not change when resuming the random clause process. This is not necessarily the case as for example some core variables may violate the self-containment property and be removed, and this may cause a chain reaction of other variables leaving the core (maybe their support is too small, or they violate the self-containment requirement). However, $\whp$ all the variables the are removed from the core when resuming the random clause process remain satellite variables, and furthermore there are very few such variables.

\begin{remark} In several papers which studied planted-solution
distributions, for example \cite{AlonKahale97,flaxman}, a similar
notion of a core appears (without the notion of satellite variables), and an analysis of the structure of the instance ($k$-colorable graph or $k$-CNF formula) induced on the non-core variables is also given. The main difference from our setting is the fact that the planted distribution is a product space, and therefore it was possible to prove that the core variables are distributed similarly to a uniformly random set of variables. In our case establishing such a property is a more challenging task. As it turns out, the approach that we take -- defining the satellite variables -- simplifies considerably the proof of this property.
\end{remark}

\subsection{The Majority Vote}\label{sec:Majority}
Given a 3CNF formula $\F$ and a variable $x$ we let $N^+(x)$ be the
set of clauses in $\F$ in which $x$ appears positively (namely, as
the literal $x$), and $N^{-}(x)$ be the set of clauses in which $x$
appears negatively (that is, as $\bar{x}$). The Majority Vote
assignment over $\F$, which we denote by MAJ, assigns every $x$
according to the sign of $|N^+(x)|-|N^-(x)|$ (TRUE if the difference
is positive and FALSE otherwise).

\begin{proposition}\label{prop:MajVoteSuccRate} Let $\F^*$ be distributed according to $\RandDistSAT$ with
$n^2p\geq d$, $d$ a sufficiently large constant. Then $\whp$ every
satisfying assignments of $\F^*$ differs from MAJ on at most
$\SmallFrac n$ variables.
\end{proposition}
\begin{Proof}
Consider the following two-step procedure to generate $F$: in the first step go over the $M=8\binom{n}{3}$ clauses and toss a coin with success probability $p_1$. We take the clauses that were chosen
and put them first, ordered at random. Call $F_1$ this first part (and respectively define $F_1^*$ in our standard way, i.e., by scanning sequentially the clauses of $F_1$  and including those whose addition leaves the formula satisfiable.). Observe that $F_1$ is distributed according to ${\cal{P}}_{n,p_1}$. Then in the second round, every clause that was not chosen in the first round is included with probability $p_2$, and the chosen clauses are ordered at random
and then concatenated after $F_1$. Call $F_2$ this last part.
At the end of this subsection we prove that $F=F_1 \cup F_2$ is distributed according to $\RandDist$
when $p=p_1 +(1-p_1)p_2$. Therefore we may think of $F$ as generated in two steps (with the suitable choice of $p_1,p_2$). We will use this technique to prove several other properties as well.

Let $d_0$ be the constant promised in Proposition \ref{prop:CoreSize}, and
choose $d \geq 200d_0$. Set $p_1 = p/200$. By the choice of $d_0$ and Proposition \ref{prop:CoreSize} $\whp$ all but $\SmallFrac n$ variables are frozen in $F_1^*$, and w.l.o.g assume that they all take the value TRUE. Further observe that $\whp$ at this point all but $\SmallFrac n$ variables appear in no more than say $n^2p/30$ clauses (in $F_1$ distributed according to ${{\cal{P}}_{n,p_1}}$, and therefore also in $F_1^*$). This is because every variable $x$ is expected to appear in $F_1$ in $p_1\cdot 8\binom{n}{2}\leq 4n^2p_1 = n^2p/50$ clauses. These appearances are independent (binomially distributed), therefore one can apply the Chernoff bound for example to bound the probability that $x$ appears in more than $n^2p/30$ clauses, which will be $\SmallFrac$. This in turn gives that the expected number of such variables is $\SmallFrac n$. To obtain concentration around this value, consider an ordering on the $M$ clauses and let $X_i$ be an indicator random variable which is 1 iff clause $i$ appeared in the first round. Let $f(X_1,X_2,\ldots,X_M)$ be a function which counts the number of variables that appear in more than  $n^2p/30$ clauses in $F$. As claimed, $E[f]=\SmallFrac n$, and $f$ satisfies the Lipschitz condition with difference 3: for every $i$ and every two assignments $a=(a_1,\ldots,a_{i-1},a_i,a_{i+1},\ldots,a_M)$ and $a'=(a_1,\ldots,a_{i-1},a'_i,a_{i+1},\ldots,a_M)$ of values to $X_1,\ldots,X_M$ (that possibly differ on the $i^{th}$ coordinate), it holds that $|f(a)-f(a')| \leq 3$ (every clause contains three variables). Using the method of bounded differences (e.g., Theorem 7.4.3 of \cite{TheProbMethod}) it follows that $f$ is concentrated around its expected value.

Let $Z$ be then the set of frozen variables that appear in at most
$n^2p/30$ clauses of $\F_1$. Recall that we have assumed w.l.o.g. that
they all froze to TRUE. By the above discussion together with Proposition
\ref{prop:CoreSize} $\whp$ $$|Z|\geq (1-\SmallFrac)n-\SmallFrac n
\geq 0.999n.$$ Now let us consider the second iteration of coin
flips. Fix $x \in Z$, observe that every clause containing $x$
positively, if chosen in the second round will be included in
$F^*$. There are at least $4\binom{|Z|-1}{2}- n^2p/30$ such clauses
with the other two variables from $Z$ -- call them ``good" clauses. As for clauses where $x$ appears negatively,
and the other two variables are in $Z$, there are only at most
$3\binom{|Z|-1}{2}$ clauses such that if chosen will be included
(since one way of negating the variables in $Z$  results in a FALSE clause
on frozen variables) -- call them ``bad" clauses. In addition there are
at most $8(n-|Z|)n$ clauses, containing $x$ and at least one variable outside $Z$,
that we don't say anything about, but let us
adversarially assume that $x$ appears in all of them negatively, and if
chosen are included in $F^*$ (they are also part of the bad clauses).

In expectation, $p_2\cdot \left(4\binom{|Z|-1}{2}-n^2p/30\right)
\geq 1.8n^2p$ good clauses containing $x$ will be chosen in the second
round, and
$p_2\cdot \left(3\binom{|Z|-1}{2}+8(n-|Z|)n\right) \leq 1.6n^2p$
bad clauses. (Recall that $199p/200 \le p_2\le p$.)

Suppose that in the $n^2p/30$ clauses from the first round also $x$ appears
negatively. To conclude, for the majority vote of $x$ to be wrong it must
have been the case that the number of good clauses containing $x$ or the
number of bad clauses containing $x$ deviates by
at least $(1.8-1.6-1/30)n^2p/2$ from its expectation. But since both are
binomially distributed with expectation $\Theta(n^2p)$, this happens with
probability $\SmallFrac$. Using the linearity of expectation all but
$\SmallFrac n$ of the variables in $Z$ are expected to have a ``proper"
gap. To obtain concentration around this
value we use again the method of bounded differences, similarly to
hat has been used earlier in the proof.
Finally observe that $|Z|\geq (1-\SmallFrac)n$, and therefore $|Z|-\SmallFrac n =  (1-\SmallFrac)n$ as required.
\end{Proof}
\subsubsection{Justifying the two-step distribution.}
Let ${{\cal{P}}_{n,p_1,p_2}}$ be the distribution of the two-step process. For brevity, set
${\cal P}_1=\RandDist$, ${\cal P}_2={{\cal{P}}_{n,p_1,p_2}}$. Let us now prove that
${\cal P}_1$ and ${\cal P}_2$ are identical for
$p=p_1+(1-p_1)p_2$. Let $\sigma$ be an ordered list of $|\sigma|=r$ clauses. Then
$$
Pr_{{\cal P}_1} [\mbox{ get list $\sigma$}]=
\frac{p^{r}(1-p)^{M-r}}
     {r!}.
$$
On the other hand,
\begin{eqnarray*}
Pr_{{\cal P}_2} [\mbox{ get list $\sigma$}] &=&
\sum_{i=0}^r \frac{p_1^i(1-p_1)^{M-i}}{i!}\,\cdot\,
             \frac{p_2^{r-i}(1-p_2)^{M-r}}{(r-i)!}\\
&=& (1-p_1)^M p_2^r (1-p_2)^{M-r}\sum_{i=0}^r
\frac{\left(\frac{p_1}{1-p_1}\right)^i}{i!}\,\cdot\,
\frac{\left(\frac{1}{p_2}\right)^i}{(r-i)!}\\
&=& (1-p_1)^M p_2^r (1-p_2)^{M-r} \frac{1}{r!}\sum_{i=0}^r
{r\choose i} \left(\frac{p_1}{(1-p_1)p_2}\right)^i\\
&=& (1-p_1)^M p_2^r (1-p_2)^{M-r} \frac{1}{r!}
\left(1+\frac{p_1}{(1-p_1)p_2}\right)^r\\
&=& \frac{((1-p_1)(1-p_2))^{M-r}(p_1+p_2-p_1p_2)^r}{r!}\,.
\end{eqnarray*}

Choosing $p_1,p_2$ to satisfy $p_1+p_2-p_1p_2=p$, we conclude
that the distributions ${\cal P}_1$ and ${\cal P}_2$ are
indeed identical.

\subsection{Proof of Theorem \ref{thm:StructOfSolutionSpace}}
Theorem \ref{thm:StructOfSolutionSpace} follows from Proposition \ref{prop:SizeOfConnectedComp}
which implies that all but $\SmallFrac n$ of the variables are frozen.
Therefore, there are at most $\exp\{{\SmallFrac n}\}$ possible ways to
set the assignment of the remaining variables. Furthermore, Proposition \ref{prop:SizeOfConnectedComp}
describes the formula induced by the non-frozen variables.

\section{Proof of Theorem \ref{thm:PlytimeAlg}}\label{sec:ProofOfThm}
\begin{figure*}[!htp]
\begin{center}
\fbox{
\begin{minipage}{\textwidth}{\textbf{\textsf{SAT}$(\F,t)$}}\\
\noindent\underline{\textsl{Step 1: \textsf{Majority Vote}}}

1. $\pi_1\leftarrow$ Majority Vote over $\F$.

\noindent\underline{\textsl{Step 2: Reassignment}}

2. \textbf{for} $i=1$ \texttt{to} $\log n$

\noindent3.  \text{ }\text{ }\text{ }\textbf{for all}
$x\in V$

\noindent 4.\text{ }\text{ }\text{ }\text{ } \text{ }\text{ }\textbf{if} $x$
supports less than $2t/3$ clauses w.r.t. $\pi_{i}$ \textbf{then} $\pi_{i+1}\leftarrow \pi_{i}$ with $x$ flipped.

\noindent 5. \text{ }\text{ }\text{ }\textbf{end for.}

\noindent 6. \textbf{end for.}

\noindent\underline{\textsl{Step 3: Unassignment}}

\noindent 7. set $\psi_1=\pi_{\log n}$, $i=1$.

\noindent 8. \textbf{while} $\exists x$ s.t. $x$ supports less than
$t$ clauses w.r.t. $\psi_i$

\noindent 9.\text{ }\text{ }set $\psi_{i+1}
\leftarrow \psi_i$ with $x$ unassigned.

\noindent 10.\text{ }\text{ }\text{ }$i\leftarrow i+1.$

\noindent 11. \textbf{end while.}

\noindent\underline{\textsl{Step 4. Unit Clause Propagation}}

\noindent 12. Let $\xi$ be the final partial assignment obtained
at Step 3.

\noindent 13. Remove all clauses which are satisfied by $\xi$, and all FALSE-literals from the remaining clauses.

\noindent 14. Run the unit-clause-propagation algorithm on the resulting instance.

\noindent\underline{\textsl{Step 5: Exhaustive Search}}

\noindent 15. Let $F'$ be the formula remaining after the
unit-clause-propagation of Step 4 terminates.

\noindent 16. Exhaustively search and satisfy $\F'{out}(A,\xi)$, component by component.
\end{minipage}
}
\end{center}\caption{The algorithm \textsf{SAT}} \label{fig:SATAlg}
\end{figure*}
In this section we prove that the algorithm \textsf{SAT}, which is
described in Figure \ref{fig:SATAlg}, meets the requirements of
Theorem \ref{thm:PlytimeAlg}. The main principles underlying \textsf{SAT} were designed with
the planted distribution in mind (see \cite{flaxman} for example). An additional ingredient that we add is a unit-clause-propagation step. Given a 1-2-3-CNF formula (namely a formula which contains clauses of size 1,2 and 3), the unit-clause-propagation is the following simple heuristic:

\begin{quote}
\it while there exists a clause of size 1, set the variable appearing in
this clause in a satisfying manner, remove this clause and all
other clauses satisfied by this assignment, and remove the FALSE literals
of the variable from other clauses.
\end{quote}

We say that $\F^*$ is \emph{typical} in $\RandDistSAT$ if
Propositions
\ref{prop:SizeOfConnectedComp} and  \ref{prop:MajVoteSuccRate} hold. The
discussion in Section \ref{sec:PropOfRandomInst} guarantees that
indeed $\whp$ $\F^*$ is typical. Therefore, to prove Theorem
\ref{thm:PlytimeAlg} it suffices to consider a typical $\F^*$ and
prove that \textsf{SAT} (always) finds a satisfying assignment for
$\F^*$. As the parameter $t$ for \textsf{SAT} we use the $t$
promised in Proposition \ref{prop:SizeOfConnectedComp}.

We let $\Core$ be the $t$-core promised in
Proposition~\ref{prop:SizeOfConnectedComp}, $\Sat$ its satellite variables, and $\varphi$ be the
satisfying assignment w.r.t. which $\Core$ is defined. In all the
following propositions we assume $\F^*$ is typical (we don't
explicitly state it every time for the sake of brevity).

\begin{proposition}\label{prop:ReassigmentCorrect} Let
$\psi_1$ be the assignment defined in line 7 of \textsf{SAT}.
Then $\psi_1$ agrees with $\varphi$ on the assignment of all
variables in $\Core$.
\end{proposition}
\begin{Proof}
Let $B_i$ be the set of core variables whose assignment in $\pi_i$
disagrees with $\varphi$ at the beginning of the $i^{th}$ iteration
of the main for-loop -- line 2 in \textsf{SAT}. It suffices to prove
that $|B_{i+1}|\leq |B_i|/2$ (if this is true, then after $\log n$
iterations $B_{\log n}=\emptyset$). Observe that by Proposition
\ref{prop:MajVoteSuccRate}, $|B_0|\leq n/10^7$ (as the Majority Vote error-rate
$\SmallFrac$ can be made arbitrarily small). By contradiction,
assume that not in every iteration $|B_{i+1}|\leq |B_i|/2$, and let
$j$ be the first iteration violating this inequality. Consider a variable $x\in B_{j+1}$. If also
$x \in B_j$, this means that $x$'s assignment was not flipped in the
$j^{th}$ iteration, and therefore, $x$ supports at least $2t/3$
clauses w.r.t. $\pi_j$. Since $\Core$ is $t/3$-self-contained, at
least $2t/3-t/3 = t/3$ of these clauses contain only core
variables. Since the literal of $x$ is true in all these clauses,
but in fact should be false under $\varphi$, each such clause must
contain another variable on which $\varphi$ and $\pi_j$ disagree,
that is another variable from $B_j$. If $x\notin B_j$, this means
that $x$'s assignment was flipped in the $j^{th}$ iteration. This is
because $x$ supports less than $2t/3$ clauses w.r.t. $\pi_j$. Since
$x$ supports at least $t$ clauses w.r.t. $\varphi$ ($t$-expanding
property of the core), it must be that in at least $t-2t/3=t/3$ of
them, the literal of some other core variable evaluates to TRUE
(not FALSE as it should be in $\varphi$). Letting $U=B_j \cup B_{j+1}$, there are at least
$t/3\cdot|B_{j+1}|$ clauses containing at least two variables from
$U$ (every clause is counted exactly once as the supporter of a clause
is unique). Using our assumption, $|B_{j+1}|\geq |B_j|/2$, we obtain $|U|=|B_j \cup B_{j+1}|\leq |B_j|+|B_{j+1}|\leq
3|B_{j+1}|$, therefore $t/3\cdot|B_{j+1}| \geq (|U|/3)\cdot t/3 = (t/9)|U|$.
Finally,
\begin{itemize}
  \item $|B_j|\leq n/10^7$ (because $B_0$ is already small enough, and by our assumption the sets $B_1,B_2,\ldots B_j$ only decrease in size),
  \item $|B_{j+1}|$ may exceed $n/10^7$, in which case we consider w.l.o.g. only the first $n/10^7$ variables (this is in line with our assumption $|B_{j+1}|\geq |B_j|/2$),
  \item $|U|\leq 3|B_{j+1}|\leq 3n/10^7\leq n/10^6$,
  \item there are $t|U|/9$ clauses containing two variables from $U$.
\end{itemize}
The last two items contradict the $t/10$-proportionality of $\F^*$.

\end{Proof}

\begin{proposition}\label{prop:UnassigmentCorrect} Let
$\xi$ be the partial assignment defined in line 12 of \textsf{SAT}.
Then all assigned variables in $\xi$ are assigned according to
$\varphi$, and all the variables in $\Core$ are assigned.
\end{proposition}

\begin{Proof}
By Proposition \ref{prop:ReassigmentCorrect}, $\psi_1$ coincides
with $\varphi$ (the satisfying assignment w.r.t. which
$\Core$ is defined) on $\Core$. Furthermore, by the definition of
$t$-core, every core variable supports at least $t$ clauses w.r.t.
$\varphi$, and also w.r.t. $\psi_1$ (the assignment at hand before
the unassignment step begins). Hence all core variables survive the
first round of unassignment. By induction it follows that the core
variables survive all rounds. Now suppose by contradiction that not
all assigned variables are assigned according to $\varphi$ when the
unassignment step ends. Let $U$ be the set of variables that remain
assigned when the unassignment step ends, and whose assignment
disagrees with $\varphi$. Every $x\in U$ supports at least $t$
clauses w.r.t. to $\xi$ (the partial assignment defined in line 12 of
\textsf{SAT}), but each such clause must contain another variable on
which $\xi$ and $\varphi$ disagree (since $\varphi$ satisfies this clause). Thus, we have $t\cdot |U|$ clauses each
containing at least two variables from $U$ (again no clause is counted twice as the support of a clause is unique). Since $U\cap \Core
=\emptyset$ (by the first part of this argument) and $|\Core|\geq
(1-\SmallFrac)n$ it follows that $|U|\leq \SmallFrac n <
n/10^6$, contradicting the $t/10$-proportionality of
$\F^*$.
\end{Proof}

\begin{proposition}\label{prop:UnitClauseCorrect} By the end of the unit-clause propagation step
all the variables which get assigned are assigned according to
$\varphi$, furthermore the set of satellite variables $\Sat$ is assigned.
\end{proposition}
\begin{Proof}
The proof is by induction on the iterations of the unit clause propagation.
The base case are clauses of the form $(x \vee \ell_z \vee \ell_y)$ where $\ell_z,\ell_y$ are FALSE literals under $\xi$ and $x$ is unassigned. By the previous proposition, $\xi$ can be extended to a satisfying assignment of $F$, but every such extension must set $x=TRUE$. This is exactly what the unit clause propagation does. The step of the induction is proven similarly to the base case.

Now to the satellite variables. The previous proposition gives that $\Core$
remains assigned according to $\varphi$. By the definition of satellite
variables, $\Sat$ will be set in the unit clause propagation (the $i$-satellite variables will be set in iteration $i$ of the unit-clause propagation).
\end{Proof}

\begin{proposition}\label{prop:ExhastiveSearchCorrect} The exhaustive
search, Step 5 of \textsf{SAT}, completes in polynomial time with
a satisfying assignment of $\F^*$.
\end{proposition}
\begin{Proof}
By Proposition \ref{prop:UnitClauseCorrect}, the partial assignment
at the beginning of the exhaustive search step is partial to the
satisfying assignment $\varphi$ of the entire formula. Therefore the
exhaustive search will succeed. Further observe that the unassigned
variables are outside of $\Core \cup \Sat$. Proposition \ref{prop:SizeOfConnectedComp}
then guarantees that the running time of the exhaustive search will be at most polynomial.
\end{Proof}\\
Theorem \ref{thm:PlytimeAlg} follows.

\section{Proof of Proposition
\ref{prop:ExpanderSize}}\label{sec:ExpndrProof} Let $F$ be the
random $\RandDist$ instance, $F^*$ be its satisfiable part. We divide the
process of generating $F$ into two steps like in the proof of Proposition \ref{prop:MajVoteSuccRate}:
in the first round go over the $M=8\binom{n}{3}$ clauses and toss a coin
with success probability $p_1=p/2$. Take the clauses that were chosen and
put them first ordered at random. In the second round, every clause that
was not chosen, is included with probability $p_2$, $p_2$ satisfies
$p_1+(1-p_1)p_2=p$; then the included clauses are ordered at
random and concatenated after the first part. Observe that this
distribution is identical to $\RandDistSAT$ as explained before.

Let $t$ be such that $F$ (and hence also $F_1$) is $\whp$ $t$-proportional (we can choose $t=n^2p/5500$
as asserted in Proposition \ref{prop:NoDenseSubgraphs}). Also take $n^2p$ sufficiently
large so that $\F^*$ is $\whp$ $n/10^6$-concentrated (as required by Proposition \ref{prop:ExpanderSize}, and as promised to be the case $\whp$ by Proposition \ref{prop:Concentration}).

Fix $\psi$ to be some assignment (not necessarily a satisfying assignment of $F^*$), and
let $B_\psi$ be a random variable counting the number of variables
whose support in $F_1$
w.r.t. $\psi$ is smaller than $502t$. A bound of the sort $Pr[B_\psi >
n/10^7]=o(2^{-n})$ would be very useful as we can then take the union bound
over all possible assignments $\psi$. Fix some variable $x$, and
w.l.o.g. assume $x$ is TRUE in $\psi$. There are $\binom{n-1}{2}$ clauses that $x$ supports w.r.t. $\psi$, each included w.p. $p_1$. Therefore in
expectation $x$ supports at least $n^2p_1/3= n^2p/6$ clauses.  Since the support of $x$ is distributed binomially, the probability that $x$ supports less than $t$ clauses in $F_1$ w.r.t. $\psi$ is at most $e^{-n^2p/50}$ (say, use the Chernoff bound). Finally observe that the set of clauses that $x$ supports is disjoint from the set of clauses that $y \neq x$ supports. Therefore, the probability that there are at least $n/10^7$ such variables is at most ${n\choose {n/10^7}}e^{-(n^2p/50) \cdot (n/10^7)} < 3^{-n}$ for sufficiently large $n^2p$.

In particular, $\whp$ every $\psi$ that satisfies $\F_1^*$ has the desired property. Let now $\psi$ be a satisfying assignment of $\F_1^*$ such that $B_\psi\leq n/10^7$, and consider the
following procedure which, as we shall prove, produces a large
$500t$-expanding set $Z$ in $\F_1^*$ (and therefore also in $\F^*$ which
contains $\F_1^*$). When using the notation $\F[A]$ for a formula $\F$ and
a set of variables $A$ we mean all clauses in $\F$ in which all three variables belong to $A$.
\begin{figure}[!htp]
\begin{center}
\fbox{
\begin{minipage}{\textwidth}\it
\begin{enumerate}
\item set $Z_{0} = V \setminus \{x\in V: x \text{ supports less
than $502t$ clauses in $\F_1^*$ w.r.t. $\psi$}\}$; $i=0$.
\item \textbf{while} there exists a variable $a_{i} \in Z_{i}$
that supports less than $500t$ clauses in $\F_1[Z_{i}]$
\textbf{do} $Z_{i+1}=Z_{i} \setminus \{a_{i}\}$; $i\leftarrow i+1$.
\item let $a_{r}$ be the last variable removed in step 2. Define
$Z=Z_{r+1}$.
\end{enumerate}\rm
\end{minipage}
}
\end{center}\caption{Building a $t$-expanding set} \label{fig:ExpandingSet}
\end{figure}

Clearly, $Z$ is $500t$-expanding in $\F_1$ (by the construction). It
remains to prove that $Z$ is large. By our assumption on
$B_\psi$ step 1 removes at most $n/10^7$ variables, let $A$ be
those variables. It remains to prove that in the iterative step not
too many variables were removed. Suppose by contradiction that in
the iterative step more than $n/10^7$ variables were removed, and
consider iteration $j=n/10^7$ and the set $W=\{a_1,\ldots,a_{j}\}$ ($a_i \in W$ is defined in line 2 of Figure \ref{fig:ExpandingSet}). Every $a_i \in W$ appears in more than
$502t-500t=2t$ clauses in which at least another variable belongs to
$U=W\cup A$ (by the choice of $Z_0$ and the condition in line 2 that caused $a_i$ to be removed). Therefore, by iteration
$n/10^7$, the set $U$ contains at most $n/10^7+n/10^7
\leq n/10^6$ variables, and there are more than $2t \cdot |W| \geq 2t \cdot |U|/2= t|U|$ clauses
containing at least two variables from $U$ (no clause is counted twice as the support
of a clause is unique). This contradicts the
$t$-proportionality of $F_1$. To conclude, $|Z|\geq
\left(1-10^{-6}\right)n\geq 0.99n$ as required. Observe that $|W| \geq |U|/2$ by our assumption on the size
of $A$ and by the choice of $j=n/10^7$.

\medskip

It follows that $\whp$ for every satisfying assignment $\psi$ of
$F_1^*$ there exists a $500t$-expanding set $Z$ of variables of
cardinality $|Z|\ge 0.99n$. W.l.o.g. we can take $Z$ to be maximal
such set.

Observe that $Z$ and $F_1^*$ satisfy the conditions of Proposition
\ref{prop:UniqueOfSat} (that is, $\psi$ is a satisfying
assignment, $F_1^*$ is $t$-proportional and $n/10^6$-concentrated) and
therefore $Z$ is frozen in $F_1^*$; w.l.o.g. assume that all variables in
$Z$ froze to TRUE.
Since all variables of $Z$ are frozen in $F_1^*$, we can take {\em
the same} $Z$ for every satisfying assignment $\psi$ of $F_1^*$.

So let $Z$ be as above, $|Z|\ge 0.99n$.  Now we
consider the second round of coin tosses, call the chosen clauses $F_2$.
We prove that after adding them, with probability $1-o(2^{-n})$ $Z$ extends
to a $t$-expanding set
$Z'$, $Z\subseteq Z'$, of the required size ($|Z'|\geq (1-\SmallFrac)n$).
Fix some variable $x \notin Z$ and observe that
$x$ supports $\binom{|Z|}{2}$ clauses, where $x$ appears without
negation and the other two variables are in
$Z$ and appear as negated. Since $x \notin Z$, we know
that in the first iteration at most $500t$ such clauses were included.
In expectation, $F_2$ contains at least $p_2\left(
\binom{|Z|}{2}-500t\right) \geq n^2p/5 \geq 1000t$ such clauses (this is due to $p_2\ge p/2$).
If indeed at least $500t$ clauses are included then $Z \cup \{x\}$
is a $500t$-expanding set. The probability that less than $500t$
of them were included is $\SmallFrac$ (again, Chernoff bound).
We can argue similarly about the number of clauses in $F_2$,
containing $x$ and two variables from $Z$, where all three
variables appear as negated.

Call a variable $x$ {\em good} if it participates it at least
$500t$ clauses in $F_2$ where the other two variables are from $Z$
and are negated and $x$ is not negated, and also in at least
$500t$ clauses in $F_2$ where the other two variables from $Z$ and
all three variables are negated; otherwise $x$ is called {\em
bad}. Observe that for every good $x$, for every satisfying
assignment $\psi$ of $F_1^*\cup F_2^*$, we can add $x$ to $Z$,
regardless of whether $\psi$ sets $x$ to TRUE or FALSE. The
above argument shows that the expected number of bad variables is
$\SmallFrac n$. Applying standard concentration techniques, we can
derive that $\whp$ the number of bad variables is $\whp$
$\SmallFrac n$ as well.

And so we have proven that $\whp$ there exists a $500t$-extending set $Z'$ (which contains $Z$) and
$|Z'|=|Z|+(1-\SmallFrac)|V\setminus Z| \geq (1-\SmallFrac)|V|$.

\medskip

For conclusion, we have shown that there exists a $500t$-expanding set
$Z'$ in $\F^*$ of cardinality $|Z'|=(1-\SmallFrac)n$ w.r.t. $\psi$,
where $\psi$ is some satisfying assignment of $\F^*$ (in fact this
is true w.r.t. all satisfying assignments of $F^*$ by the frozenness
property). Scaling everything down (setting $t'=500t$), $Z'$ is
$t'$-expanding and (at least) $t'/500$-proportional. This
completes the proof of the proposition.

\begin{remark}\label{rem:t/500-prop}
Note that here we proved $t'/500$-proportionality, which is stronger than what we are required to prove
($t'/10$-proportionality). In general, we prefer clear and shorter
presentation over optimizing the constants in the proofs. Later we
will use this slackness in other proofs that rely on this one.
\end{remark}

\subsection{Proof of Proposition
\ref{prop:CoreSize}}\label{sec:CoreSizeProof} Let $Z$ be the
$t$-expanding set promised by Proposition \ref{prop:ExpanderSize}.
Consider the procedure in Figure \ref{fig:Core}, which shall produce a
$t'$-core (for $t'=10t/11$). Recall that using the notation $\F[A]$ for formula $\F$ and set of variables $A$ we mean all clauses in $\F$ in which all three variables belong to $A$.
\begin{figure*}[!htp]
\begin{center}
\fbox{
\begin{minipage}{\textwidth}\it
\begin{enumerate}
\item set $H_{0} = Z$ and $i=0$.
\item \textbf{while} there exists a variable $a_{i} \in H_{i}$ that:
 \begin{itemize}
   \item $a_i$ appears in more than $t/11$ clauses in $\F\setminus \F[H_{i}]$, \textbf{or},
   \item $a_i$ supports less than $10t/11$ clauses in $\F[H_{i}]$,
 \end{itemize}
\textbf{do} $H_{i+1}=H_{i} \setminus \{a_{i}\}$.
\item let $a_{r}$ be the last variable removed in step 2. Define
$\Core=H_{r+1}$.
\end{enumerate}\rm
\end{minipage}
}
\end{center}\caption{Building a $t$-core} \label{fig:Core}
\end{figure*}

First let us explain why indeed $\Core$ is a $t'$-core. By its construction $\Core$
is $10t/11$-expanding (or $t'$-expanding). Further, $\Core$ is $t/11$-self-contained, or $t'/10$-self-contained
(which also implies $t'/3$-self-contained as $1/3>1/10$).

\begin{remark}\label{rem:t/6-self-contained}
By the definition of a core we are required to prove only $t'/3$-self-contained, but we shall need this slackness in the proof of Proposition \ref{prop:SizeOfConnectedComp}.
\end{remark}

It remains prove that $|\Core|\geq(1-\SmallFrac)n$. By
Proposition \ref{prop:ExpanderSize}, $|H_0| \geq (1-e^{-n^2p/c_1})n$
for some constant $c_1>0$ independent of $n,p$. Let $A=V\setminus H_0$, and note that $|A| \leq e^{-n^2p/c_1}n$. Suppose that the iterative procedure (line 2) removed more than $e^{-n^2p/c_1}n$ variables.
Consider iteration $j=e^{-n^2p/c_1}n$ and the set $W=\{a_1,\ldots,a_j\}$ ($a_i \in W$ is defined in line 2 of Figure \ref{fig:Core}). Define $U = W \cup A$. One possibility for the removal of $a_i$ is that it appears in
at least $t/11$ clauses in which at least another variable belongs to
$U$. Another is that $a_i$ supports less than
$10t/11$ clauses w.r.t. $\F[H_{i}]$. In the latter case $a_i$ must support at
least $t-10t/11=t/11$ clauses in $\F\setminus \F[H_{i}]$ (by the choice of $a_i\in Z$).
In any case $a_i$ appears in at least $t/11$ clauses with at least another variable from $U$.
Therefore, by  iteration $j$, there exists a set
$U$ containing at most $2e^{-n^2p/c_1}n<< n/10^6$ variables, and there are
at least $(t/33) \cdot |W| \geq (t/33) \cdot (|U|/2) = t|U|/66$ clauses containing at least two
variables from it (we divide $t/11$ by 3 as a clause could have been counted three times). This however contradicts the $t/500$-proportionality of $F$ (recall that in Proposition \ref{prop:ExpanderSize}, when proving the existence of a $t$-expanding set $Z$, we actually proved that $\F$ is
$t/500$-proportional -- Remark \ref{rem:t/500-prop}).

Finally, the core variables are frozen as they are a subset of $Z$, and $Z$
-- the $t$-expanding set -- is frozen.

\section{Proof of Proposition \ref{prop:SizeOfConnectedComp}}\label{sec:ProofOfPropSizeOfConnectedComp}
Let $F$ be the random $\RandDist$ instance, $F^*$ its satisfiable part. We divide the
process of generating $F$ into two steps like in the proof of Proposition \ref{prop:MajVoteSuccRate}:
in the first round go over the $M=8\binom{n}{3}$ clauses and toss a coin with success probability $p_1=p/2$. Take the clauses that were chosen and put them first. In the second round, every clause that was not chosen is included with probability $p_2$, where $p_2$ satisfies $p_1+(1-p_1)p_2=p$. Let $F_1^*$ be the part of $F^*$ that corresponds to the first iteration ($F_1^*$ is distributed according to ${{\cal{P}}^{{\rm sat}}_{n,p_1}}$). Let $F_2$ be the clauses that were chosen in the second round.

By Proposition \ref{prop:CoreSize}, if we take $n^2p$ to be sufficiently large, then $F_1^*$ has $\whp$ a $t$-core $\Core$ w.r.t. to a satisfying assignment $\psi$ with the following properties (the last property did not appear in Proposition \ref{prop:CoreSize}, we define and justify it immediately after):
\begin{itemize}
  \item $\F_1$ is $t/500$-proportional (Remark \ref{rem:t/500-prop}).
  \item $\Core$ is $t/10$-self-contained and not only $t/3$-self-contained (Remark \ref{rem:t/6-self-contained}).
  \item $|\Core| \geq (1-\SmallFrac)n$.
  \item $\Core$ is frozen.
  \item $F_1^*$ is bounded.
\end{itemize}
We say that a formula $F$ is \emph{bounded} if no variable appears in more than $n$ clauses. In $F_1$ every variable is expected to appear in $O(n^2p)$ clauses, and we may assume that $n^2p = O(n^{1/2})$ (if not, then in particular $\whp$ $\Core=V$ and the entire discussion in this section is unnecessary). Standard calculations then show that $\whp$ no variable appears in more than $n$ clauses of $F_1$.

We now discuss what happens to $\Core$ in the second round, that is when adding $F_2$. We will be interested in  large connected components of $F_1$ whose vertices are not in $\Core$ (Proposition \ref{prop:SatteliteShatter}), and also in vertices that may leave $\Core$ due to $F_2$ (Propositions \ref{prop:Core'} and \ref{prop:IfNotCoreThenSat}). The key to understanding the transformation that $\Core$ and the connected components undergo lies in the notion of satellite variables.

First observe that $F_1$ is $\whp$ $t/500$-proportional, and therefore also is $F_2$ (as they are almost identically distributed, and there is enough slackness in the choice of constants to accommodate this difference). Hence $\whp$ $F=F_1 \cup F_2$ is $t/250$-proportional (and so is $F^*$). Assume that this is the case.

\begin{proposition}\label{prop:Core'} If after the second round $F^*$ remains $t/250$-proportional, then there exists a satisfying assignment $\psi$ of $F^*$ and a set $\Core'\subseteq \Core$ of variables which is a $t/2$-core of $F^*$ w.r.t. $\psi$. Furthermore, $|\Core'| \geq (1-\SmallFrac)n$.
\end{proposition}
\begin{Proof}
We call a variable $x \in \Core$ \emph{dirty} if in $F_2$ there exists a clause $C$ containing $x$ and some variable not in $\Core$. Let $D$ be the set of dirty variables. For a specific $x$, there are $\SmallFrac n^2$ clauses such that if chosen to $F_2$ will make $x$ dirty. The probability that any of them appears is at most $p_2 \cdot \SmallFrac n^2 = \SmallFrac$ (since $\SmallFrac$ is much smaller than $n^2p_2$ for sufficiently large $p$). Linearity of expectation gives $E[|D|]=e^{-\Omega(n^2p)}n$. Also observe that $D$ satisfies the Lipschitz condition with difference 3 (as every new clause can effect 3 new variables). Therefore also concentration is obtained. Let us assume from now on that indeed $|D|=\SmallFrac n$.

Consider $\Core$ after scanning $F_2$ (to complete $F^*$) and set $\Core_0 = \Core \setminus D, i=0$. Very similarly to the procedure in Figure \ref{fig:Core}, consider the following iterative procedure:
\begin{quote}
\it while there exists $x \in \Core_i$ s.t. $x$ supports less than $t/2$ clauses in $F[\Core_i]$ w.r.t. $\psi$, or appears in more than $t/6$ clauses where some variable belongs to $V \setminus \Core_i$, define $\Core_{i+1}=\Core_i \setminus \{x\}, i=i+1$.
\end{quote}
Set $\ell=e^{-cn^2p}n$, where $c$ is some constant satisfying $|D| \leq e^{-cn^2p}n$. Suppose that the iterative process reached iteration $\ell$, and let $W_\ell$ be the set of variables that were removed in iterations $1 \ldots \ell$, let $U=W_\ell \cup D$, and observe that $|W_\ell| \geq |U|/2$ by our choice of $c$. Take $x\in W_\ell$, if $x$ was removed in iteration $i$ because it appeared in more than $t/6$ clauses where some variable belongs to  $V \setminus \Core_i$, then since $x$ was part of $\Core$ to begin with, and $\Core$ was $t/10$-self-contained, then $x$ must appear in at least $t/6-t/10=t/15$ clauses in which at least another variable belongs to $U$. If $x$ was removed because it supports less than $t/2$ clauses in $F[\Core_i]$, then again, $x$ was part of $\Core$, and therefore it supports at least $t$ clauses in $F[\Core]$, and hence it must support (and, in particular, appear in) at least $t-t/2=t/2$ clauses in which some variable belongs to $U$. At any rate, every $x \in W_\ell$ appears in at least $t/15$ clauses in which at least another variable belongs to $U$. Finally,
\begin{itemize}
  \item there are $t|W_\ell|/15\cdot 1/3 \geq t|U|/90$ clauses containing at least two variables from $U$ (we divide by 3 as every clause might have been over--counted up to 3 times, and we use the fact that $|W_\ell| \geq |U|/2$),
  \item $|U| = |D|+|W_\ell| = \SmallFrac n < n/10^6$ (we used our estimate on $|D|$, and the fact that we look at the iterative process until iteration $\ell$, therefore $|W| \leq \ell = \SmallFrac n$).
\end{itemize}
Combining these two facts contradicts the $t/250$-proportionality of $F^*$. Therefore if we let $W$ denote the set of variables that were removed in the iterative step, in all iterations, then $\whp$ $|W| \le\ell$.
Now set $t'=t/2$, and let $\Core'=\Core \setminus \{D \cup W \}$. We have shown that the set $\Core'$ is a $t'$-core of the required size. Further, $F^*$ is (at least) $t'/10$-proportional as required by Proposition \ref{prop:SizeOfConnectedComp}.

Finally observe that $\Core$ is frozen and hence $\Core' \subseteq \Core$ is frozen too. Therefore although $\Core'$ is defined w.r.t. $\psi$, it will be a core of $F^*$ regardless of which satisfying assignments survive at the end (as it will be a core w.r.t. all $F^*$'s satisfying assignments, and at least one is guaranteed to survive).
\end{Proof}

\begin{proposition}\label{prop:IfNotCoreThenSat} Let $\Sat'$ be the satellite variables of $\Core'$. If $F^*$ is $t/10$-proportional then $\Core \setminus \Core' \subseteq \Sat'$.
\end{proposition}
\begin{Proof}
Let $A=\Core \setminus \Core'$. Let $\Sat'$ be all the satellite variables of $\Core'$, and by contradiction assume that the set $B= A \setminus \Sat'$ is non-empty.  Every $x$ in $B$ belongs to $\Core$ and therefore supports at least $t$ clauses where the other two variables appear in $\Core$. Observe that in none of these $t$ clauses the other two variables are in $\Core'\cup \Sat'$ (as otherwise $x$ is in $\Sat'$). Therefore we have found a set $B$, $|B| = \SmallFrac n \leq n/10^6$, for which there are at least $t|B|$ clauses containing two variables from $B$. This contradicts the $t/10$-proportionality of $F^*$.
\end{Proof}

In the proof of Proposition \ref{prop:SizeOfConnectedComp} we consider two ``types" of satellite variables. The first type, which we just met, are the variables in $\Core \setminus \Core'$. The second type, which we will make use of in the proof of Proposition \ref{prop:SatteliteShatter} ahead, are satellite variables of $\Core$ whose ``job" is to shatter the large connected components in the formula induced by variables not in $\Core$ (when exposing the second part of $F$). In some sense these two types represent competing processes. The one is variables leaving $\Core$, but still remaining satellite variables, the other is new variables attaching to $\Core$ as satellite variables.

\medskip

Recall our notation $\F_{out}(A,\varphi)$ ($A$ a set of variables, $\varphi$ an assignment) which stands for the subformula of $\F$ which is the outcome of the following procedure: set the variables in $A$ according to $\varphi$ and simplify $\F$ (by simplify we mean remove every clause that contains a TRUE literal, and remove FALSE literals from the other clauses). The connected components of a formula $\F$ are the sub-formulas $\F[C_1],\ldots,\F[C_k]$,
where $C_1,C_2,\ldots,C_k$ are the connected components in the
graph $G_\F$ induced by $\F$ (the vertices of $G_\F$ are the variables, and two
variables are connected by an edge if there exists some clause containing them
both).

\begin{proposition}\label{prop:SatteliteShatter} Let $\Core$ be a $t$-core of $F_1^*$, let $\Sat$ be the set of all satellite variables of $\Core$ in $F^*$, and let $\psi$ be a satisfying assignment of $F^*$. Then the largest connected component in $F^*_{out}[\Core \cup \Sat,\psi]$ is $\whp$ of
size at most $\log n$.
\end{proposition}

First let us show why Proportion \ref{prop:SatteliteShatter}
completes the proof of Proposition \ref{prop:SizeOfConnectedComp}.
Since the proposition is true for $F$ it is true, by monotonicity,
for $F^*$. We take $\Core'$ for the core to be given by Proposition
\ref{prop:SizeOfConnectedComp}, and denote by $\Sat'$ its satellite variables.
Observe that (under the assumption of proportionality) $\Core
\subseteq \Core' \cup \Sat'$, and hence by the definition of
satellite variables $\Sat \subseteq \Sat'$. In particular  $\Core \cup \Sat
\subseteq \Core' \cup \Sat'$.

\subsection{Proof of  Proposition \ref{prop:SatteliteShatter}}
Let us refine the process of generating $F$: first we generate
$F_1$ (and $F_1^*$), and fix $\Core$ according to $F_1^*$. Then in
the second round ($F_2$) first toss the coins of clauses $C$ s.t. at
most one literal in $C$ belongs to $\Core$, call $J \subseteq F_2$
the set of clauses that were chosen. Finally toss the coins of the
other clauses (the ones that were not picked in the first step and
contain at least two variables from $\Core$), call $K \subseteq
F_2$ the set of clauses that were chosen. In this new terminology
$F=F_1 \cup J \cup K$, and set $F'=F_1 \cup J$. To prove Proposition
\ref{prop:SatteliteShatter} it suffices to consider only trees of
size $\log n$ in $F'$. This is because $(a)$ every
connected component of size at least $\log n$ contains a tree of size $\log
n$, and $(b)$ only the clauses of $F'$ may contribute edges to the
connected components of $F^*_{out}[\Core \cup \Sat,\psi]$.

We will prove Proposition \ref{prop:SatteliteShatter} as follows:
fix an arbitrary tree $T$ on $r$ vertices, and let $V(T)$ denote its
set of vertices. The following two conditions are necessary for $T$
to belong to $F^*_{out}[\Core \cup \Sat,\psi]$:
\begin{itemize}
    \item $A=\{$there exists a subformula of $F'$ that induces
    $T\}$,
    \item $B=\{$the clauses in $K$ do not prevent the following from holding: $V(T)\cap \Sat=\emptyset \}$.
\end{itemize}
The probability that
$F^*_{out}[\Core \cup \Sat,\psi]$ contains a tree of size at least $r$
is at most
$$\sum_{T:|V(T)|=r}Pr[A \wedge B]=\sum_{T:|V(T)|=r}Pr[A]\cdot Pr[B|A] \leq \left(\max_{T:|V(T)|=r} Pr[B|A]\right)
\cdot \left(\sum_{T:|V(T)|=r}Pr[A]\right) \equiv q \cdot h.$$ Our
next goal is to bound $q$ and $h$, and then to show that $q \cdot
h=o(n^{-3})$ for $r = \log n$.
% (according to our new definition of
%$\whp$, see Remark \ref{rem::RedfinigWhp}).
In fact we shall prove
that $q \cdot h =o(n^{-\Omega(n^2p)})$ for $r=\log n$. The next two
lemmas establish the desired bounds (we use $d=n^2p$).

\begin{lemma}\label{lem:SatelliteLem1} $h=\sum_{T:|V(T)|=r}Pr[A]\leq n(100d)^r$.
\end{lemma}

\begin{lemma}\label{lem:SatelliteLem2} $q=\max_{T:|V(T)|=r} Pr[B|A] \leq e^{-dr/8}$.
\end{lemma}
To conclude, for $r=\log n$,
$$q \cdot h \leq n  (100d)^{\log n} \cdot n^{-d/8} \leq  n^{1+\log(100d)-d/8}=o(n^{-\Omega(d)}).$$
The last equality is true since $d/8 >>  1+\log(100d)$ for
sufficiently large $d=n^2p$. We shall now prove the two lemmas.

\medskip

\noindent{\bf Proof of Lemma \ref{lem:SatelliteLem1}.}
The quantity $h$ to be estimated is obviously the expected number of trees of size
$r=\log n$  induced by a formula $F'\subseteq F$ and is therefore at most the expected number
of such trees induced by $F$ itself. We thus estimate from above the latter quantity.

Let $T$ be a fixed tree on $r$ variables (a tree in the regular graph sense),
and let $F_T$ be a fixed collection of clauses such that each edge of $T$ is induced by
some clause of $F_T$ -- we call such $F_T$ an \emph{inducing} set of clauses.
We say that a clause set $F_T$ is \emph{minimal} w.r.t. $T$ if by
deleting a clause from $F_T$, $T$ is not induced by the new formula anymore.
By the definition of minimality, $|F_T| \leq |E(T)| = |V(T)|-1$ (as
$T$ is a tree). In our argument we shall be interested only in
$(T,F_T)$ s.t. $F_T$ is a minimal set of clauses that induces $T$.

Given a tree $T$ of size $r$, we estimate the number of ways
to extend $T$ to a minimal inducing set $F_T$. Every clause in $F_T$
can cover either one or two edges of $T$ (it cannot cover three
edges or we have a cycle in $T$). Following the argument in
\cite{flaxman}, let $N_{T,s}$ be the number of ways to pair $2s$
edges of $T$ to form $s$ clauses in $F_T$ that cover two edges.
There are 8 ways to set the polarity of variables in every clause of $F_T$ (and
there are $r-1-s$ such clauses), and at most $n^{r-1-2s}$ ways to choose the
third variable in the $r-1-2s$ clauses that cover exactly one edge.
Using this terminology, the expected number of $r$-trees induced by a
random formula $F$, generated according to $\RandDist$ with $n^2p=d$,
is at most:
\begin{align}\label{eq:UnionBound}
\sum_{r-trees}\sum_{s=0}^{r/2}N_{T,s}8^{r-1-s}n^{r-1-2s}\left(\frac{d}{n^2}\right)^{r-1-s}
\leq \sum_{r-trees}\left(\sum_{s=0}^{r/2}N_{T,s}\right)(8d)^{r}n^{1-r}\,.
\end{align}
Our next task is to obtain useful upper bounds on the sum $\sum_{s=0}^{r/2}N_{T,s}$. To this end let us
fix a degree sequence $(d_1, . . . , d_r)$ for $T$, and consider the following procedure for \emph{properly} pairing
edges. By proper we mean that every pair of edges can be covered by a 3CNF clause; for example, we cannot pair the
edges $(x_1,x_2)$ and $(x_3,x_4)$ as they result in a 4CNF clause. For each vertex, we specify
a permutation of the edges incident to that vertex. Then we iterate through
the vertices, and for each vertex, we iterate through the edges and pair up
each unpaired edge with the edge given by the permutation associated with
the current vertex (and leave the edge unpaired if the permutation sends
the edge to itself). Any pairing of edges which can be covered by clauses
can be generated this way by choosing the permutations to transpose each
pair of edges to be covered by a single clause and to leave fixed all the other
edges. Since there are $d_i!$ different permutations for vertex $i$, we have
$$\sum_{s=0}^{r/2}N_{T,s} \leq \prod_{i=1}^r d_i!.$$
A classical result by Pr\"{u}fer is that the number of $r$-trees with degree
sequence $(d_1, . . . , d_r)$ equals $\binom{r-2}{d_1-1,\ldots,d_r-1}$ (see, for example, \cite{LovaszBook}, Section 4.1, p. 33).
There are $\binom{n}{r}$ ways to choose the $r$ vertices of the tree. So (\ref{eq:UnionBound}) is at most
$$
\sum_{d_1+\ldots_+d_r=2(r-1)}\binom{n}{r}\binom{r-2}{d_1-1,\ldots,d_r-1}\left(\prod_{i=1}^r d_i!\right)(8d)^{r}n^{1-r}
\leq \sum_{d_1+\ldots_+d_r=2(r-1)}\left(\prod_{i=1}^r d_i\right)
(8d)^rn.
$$
By convexity, for $(d_1, . . . , d_r)$ with $d_1 + . . . + d_r = 2(r - 1)$, the product $\prod_{i=1}^r d_i$ is maximized
when $d_1 =\ldots = d_r$, and so $\prod_{i=1}^r d_i \leq 2^r$.
The number of ways to choose positive integers $(d_1, . . . , d_r)$ so that $d_1 + . . . + d_r = 2(r - 1)$ is
$\binom{2r-3}{r-1}$ which is less than $2^{2r}$.
Hence, the expected number of $r$-trees induced by a
random formula $F$ is at most
$n\cdot 2^{2r}\cdot 2^r\cdot (8d)^r \le n(100d)^r$.
\hfill$\blacksquare$\bigskip

\noindent{\bf Proof of Lemma \ref{lem:SatelliteLem2}.} Fix a tree $T$ in $F'=F_1 \cup J$
on $r$ vertices (recall that $J$ is the set of clauses that contain at most one variable from $\Core$), and consider the set of clauses that have at least two variables in $\Core$, which we now toss their coins (we use $K$ to denote the set of clauses that were chosen among the latter).

Assume w.l.o.g. that
the assignment $\psi$, w.r.t. which $\Core$ is defined, is the
all-TRUE assignment.
Look at a variable $x \in V(T)$. We call a clause $(x \vee \bar{z}_1
\vee \bar{z}_2)$, $(\bar{x} \vee \bar{z}_3 \vee \bar{z}_4)$, where
the $z_i$'s are some variables in $\Core$, a type 1, respectively type 2, clause.
If clauses of both types appear in $K$ then $x$ surely belongs to $\Sat$
(and therefore $V(T) \cap \Sat \neq \emptyset$). We call $x \notin \Core$ \emph{elusive} if at
least one of the two types of clauses didn't appear in $K$.

Set $\rho=1-\SmallFrac$, since
$|\Core| \geq  \rho n$ there are at least $\binom{\rho n}{2} \geq
(\rho n)^2/3$ clauses of type 1. We assume that $F_1$ is bounded and hence every variable
appears in at most $n$ clauses, therefore at most $n$ clauses of type 1 have been included in $F_1$.
An identical argument applies to clauses of type 2.
 Note also that the clauses of $J$ cannot belong to any of the types.
Therefore the probability that no clause of type 1 belongs to $K$ is at most $(1-p_2)^{(\rho n)^2/3-n} \leq e^{-d/7}$
(here we use: $d=n^2p$ is large, $p_2=(p-p_1)/(1-p_1)$, $p_1=p/2$).
The same is true by symmetry for clauses of type 2. Let
$E_x$ be the event that $x$ is elusive, and let $P_i$ be the event
that no clause of type $i$ for $x$ appeared, $i=1,2$ (namely, $E_x =
P_1 \vee P_2$).
$$Pr[E_x] = Pr[P_1 \vee P_2] \leq Pr[P_1]+Pr[P_2]\leq 2e^{-d/7}
\leq e^{-d/8}.$$

Further observe that for $x \neq y$ the events $E_x$ and $E_y$ are
independent as they involve disjoint sets of clauses (each variable
supports its own set of clauses). Recall the events $A,B$ which were defined
above. In this terminology we just upper bounded the probability of $B$ given $A$, and therefore the following is true:
$$Pr[B|A] \leq \left(e^{-d/8}\right)^r = e^{-dr/8}.$$
Since our upper bound on $Pr[B|A]$ only depends on the fact that
$|V(T)|=r$, then also
$$q=\max_{T:|V(T)|=r} Pr[B|A]\leq e^{-dr/8}.$$
\hfill$\blacksquare$\bigskip

\section{$k$-Colorability}\label{sec:k-col}
In this section we will discuss, in a high level fashion, how one
can obtain similar results to the ones we have for $k$-SAT for the
random graph process (of $k$-colorability). Before we start our
discussion let us recall the algorithm due to Alon and Kahale for
coloring $k$-colorable graphs \cite{AlonKahale97}. The first step of
the algorithm is a spectral step; specifically a $k$-coloring of the
graph (not necessarily proper) is obtained by looking at some
eigenvectors of the graph (that hopefully reflect in some sense a
proper $k$-coloring). Then, this initial $k$-coloring is refined
using a series of combinatorial steps (very similar to our Steps
2--4 in Algorithm \textsf{SAT}), until possibly a proper
$k$-coloring is reached (or the algorithm fails). The algorithm was
analyzed on graphs drawn from the planted distribution first defined
at \cite{Kuc77} (the distribution is defined by the following
procedure: partition the vertex set into $k$ color classes of size
$n/k$ each: $V_1,V_2,\ldots,V_k$; next, include every $V_i-V_j$ edge
with probability $p$). The algorithm was shown to find $\whp$ a proper
$k$-coloring of the graph when $np\geq ck^2$, $c$ some sufficiently
large constant.

\medskip

It is possible to prove that the algorithm works also for graphs
drawn from our distribution for the same edge density (maybe the
constant $c$ is different). The main challenge is to reprove the
spectral properties of the graph. The basic idea is to notice that
$\whp$ every $k$-coloring has all of its color classes of
linear size, and also to prove discrepancy properties (similar, yet
more elaborate, to Proposition \ref{prop:NoDenseSubgraphs}). Another
crucial ingredient in the proof is establishing a similar notion of a
core (Definition \ref{def:core}).

\medskip

Unfortunately, at this point we are still unable to answer a
seemingly much simpler question: how many edges will such a graph
typically contain by the end of the process? We expect the answer to be about
$\binom{k}{2}\left(\frac{n}{k}\right)^2$ -- which would corresponds to the case
where a unique final $k$-coloring is nearly balanced.

\section{Discussion}\label{sec:Discussion}
As we already mentioned, only a vanishing proportion of
$k$-CNFs with $m$ clauses over $n$ variables are satisfiable when
$m/n$ is above the threshold. In recent years, several papers
studied different distributions over satisfiable 3CNF formulas in the
above threshold regime, more precisely some sufficiently large constant factor above the
threshold. In particular, \cite{flaxman} considered the planted 3SAT distribution, and
\cite{UniformSAT} addressed the planted and uniform distributions, both
papers developing new analytical and algorithmic techniques.
Our work joins this line of research by studying a new
distribution over satisfiable 3CNF formulas, and once again
introducing new analytical ideas to face the intricacies of
$\MainDist$. Furthermore, one interesting conclusion emerges from
combining \cite{flaxman},\cite{UniformSAT} and our result. In all
three distributions the instances show basically the same uni-cluster structure
of the solution space, and the same algorithm solves them all. This
gives rise to the following question: does forcing (in some ``natural" way) the unlikely event of being satisfiable in the above threshold
regime generally result in the structure suggested by Theorem
\ref{thm:StructOfSolutionSpace} (for clause-variable ratio greater than some sufficiently large constant)? This question
has been answered positively for the planted and uniform distributions, and in this paper for the random satisfiable 3CNF process.

%%%%%%%%%%%%%%%%%%%%%%%%%%%%%%%%%%%%%%%%%%%%%%%%%%%%%%%%%%%%%%%%%%%%%%%%%%%%%%%%%%%%%%%%%%%%%%%%%%%%%%%%%%%%%%%%%%%
%%%%%%%%%%%%%%%%%%%%%%%%%%%%%%%%%%%%%%%%%%%%%%%%% BIBLIOGRAPHY %%%%%%%%%%%%%%%%%%%%%%%%%%%%%%%%%%%%%%%%%%%%%%%%%%%%
%%%%%%%%%%%%%%%%%%%%%%%%%%%%%%%%%%%%%%%%%%%%%%%%%%%%%%%%%%%%%%%%%%%%%%%%%%%%%%%%%%%%%%%%%%%%%%%%%%%%%%%%%%%%%%%%%%%

\end{document}